\documentclass{amsart}
\usepackage{amsmath,amssymb}

\newtheorem{theorem}{Theorem}[section]
\newtheorem{proposition}[theorem]{Proposition}
\newtheorem{lemma}[theorem]{Lemma}
\newtheorem{corollary}[theorem]{Corollary}

\theoremstyle{definition}
\newtheorem{definition}[theorem]{Definition}

\theoremstyle{remark}
\newtheorem{remark}[theorem]{Remark}
\numberwithin{equation}{section}

\def\dbar{d\hspace*{-0.08em}\bar{}\hspace*{0.1em}}
\def\res#1{\,\omega(#1)}
\def\eps{\varepsilon}
\def\C{\mathbb C}
\def\N{\mathbb N}
\def\R{\mathbb R}
\def\rint#1{{}^{#1}\hspace{-1.4ex}\int}
\def\trint#1{{}^{#1}\hspace{-1ex}\int}
\DeclareMathOperator{\Diff}{Diff}
\DeclareMathOperator{\Tr}{Tr}
\DeclareMathOperator{\Res}{Res}
\DeclareMathOperator{\M}{\mathcal M}

\begin{document}
\title[The Noncommutative Residue on Conic Manifolds]
{On the Noncommutative Residue and the Heat Trace Expansion on
 Conic Manifolds}
\author{Juan B. Gil}
\address{Temple University, Department of Mathematics, Philadelphia, PA 19122}
\email{gil@math.temple.edu}
\author{Paul A. Loya}
\address{MIT, Department of Mathematics, Cambridge, MA 02139}
\email{ploya@alum.mit.edu}
%\subjclass[2000]{Primary 35C20; Secondary 58J42, 58J35, 58J40}
\keywords{Heat trace asymptotics, noncommutative residue, cone operators}

\begin{abstract}
Given a cone pseudodifferential operator $P$ we give a full
asymptotic expansion as $t\to 0^+$ of the trace $\Tr Pe^{-tA}$,
where $A$ is an elliptic cone differential operator for which the
resolvent exists on a suitable region of the complex plane. Our
expansion contains $\log t$ and {\em new} $\,(\log t)^2$ terms
whose coefficients are given explicitly by means of residue
traces. Cone operators are contained in some natural algebras of
pseudodifferential operators on which unique trace functionals can
be defined. As a consequence of our explicit heat trace expansion,
we recover all these trace functionals.
\end{abstract}
\maketitle

%%%%%%%%%%%%%%%%%%%%%%%%%%%%%%%%%%%%%%%%%%%%%%%%%%%%%%%%%%%%
\section{Introduction}
%%%%%%%%%%%%%%%%%%%%%%%%%%%%%%%%%%%%%%%%%%%%%%%%%%%%%%%%%%%%
On a smooth compact manifold $M$ without boundary we may consider
the operator algebra
$\mathcal A=\Psi^{\mathbb Z}_{c\ell}(M)/\Psi^{-\infty}(M)$,
where $\Psi^{\mathbb Z}_{c\ell}(M)$ is the algebra of classical
pseudodifferential operators of integral order and
$\Psi^{-\infty}(M)$ is its ideal of smoothing elements.
Operators of order less than $-\dim M$ are of trace class, but
the $L^2$-trace cannot be extended to the whole
algebra $\mathcal A$. However, M.~Wodzicki \cite{WoI84} and
V.~Guillemin \cite{Gu85} independently introduced a new functional
$\Res:\mathcal A\to \C$ which vanishes on commutators, so defines
a trace on $\mathcal A$. This trace functional is called the {\em
noncommutative residue} or {\em Wodzicki residue}, and is unique
in the sense that any other trace on $\mathcal A$ is a constant
multiple of $\Res$. A detailed survey about this trace can be
found for instance in \cite{Ka88}.

The noncommutative residue is closely related to the zeta function
of operators and to generalized heat trace asymptotics. In fact,
$\Res(P)$ can be defined (up to a constant) as the residue at
$z=0$ of $\Tr PA^{-z}$, for some fixed invertible
pseudodifferential operator $A$; the complex power being defined
as by Seeley \cite{See67}. Alternatively, $\Res(P)$ can be defined
as the coefficient of $\log t$ in the asymptotic expansion as
$t\to 0^+$ of $\Tr Pe^{-tA}$. Using these different but equivalent
definitions, the noncommutative residue has been extended to
various algebras of pseudodifferential operators on manifolds with
and without boundaries (cf. \cite{Fo-Go-Le-Sc96}, \cite{GrSchr01},
\cite{Lesch99}) including manifolds with singularities, cf.
\cite{Me-Ni}, \cite{Sc97}. The purpose of this paper is to obtain
a generalized heat trace expansion for cone operators, and use it
to recover the corresponding noncommutative residue(s) that can be
associated to them. 

Let $M$ be a compact manifold with boundary and let $x$ denote a
fixed boundary defining function. Let $\Psi^{m}_b(M)$ denote the
space of $b$-pseudodifferential operators over $M$, cf.
Section~\ref{TraceFunctionals}. On  manifolds with conic
singularities, spaces of the form $x^{-p}\Psi^{m}_b(M)$, $p\in\R$,
arise as natural spaces of pseudodifferential operators containing
the cone differential operators.

As originally done by R. Melrose and V. Nistor in \cite{Me-Ni} for
cusp operators, it is helpful to consider the following algebra
\[ x^{-\mathbb{Z}}\Psi^{\mathbb{Z}}_b(M) =
  \bigcup_{p\in \mathbb{Z}} \bigcup_{m \in \mathbb{Z}} x^{-p}
  \Psi^{m}_b(M). \]
In this algebra there is an ideal of smoothing operators
\[ \mathcal{I}:= x^{\infty}\Psi^{-\infty}_b(M) =
   \bigcap_{p\in \mathbb{Z}} x^{-p} \Psi^{-\infty}_b(M),\]
and we can consider the following quotient algebras:
\begin{equation*}
 \mathcal{I}_\sigma:=
 x^{\infty}\Psi^{\mathbb{Z}}_b(M)/\mathcal{I},
 \quad \mathcal{I}_\partial :=
 x^{-\mathbb{Z}}\Psi^{-\infty}_b(M)/\mathcal{I},
\end{equation*}
\begin{equation*}
\begin{split}
 \mathcal{A}_\sigma &:=
 x^{-\mathbb{Z}}\Psi^{\mathbb{Z}}_b(M)/
 x^{-\mathbb{Z}}\Psi^{-\infty}_b(M),\\
 \mathcal{A}_\partial &:=
 x^{-\mathbb{Z}}\Psi^{\mathbb{Z}}_b(M)/
 x^{\infty}\Psi^{\mathbb{Z}}_b(M),
\end{split}
\end{equation*}
\begin{equation*}
 \mathcal{A}_{\partial,\sigma} :=
 x^{-\mathbb{Z}}\Psi^{\mathbb{Z}}_b(M)/
 \{ x^{-\mathbb{Z}}\Psi^{-\infty}_b(M) +
    x^{\infty}\Psi^{\mathbb{Z}}_b(M)\}.
\end{equation*}
These quotients `separate' the filtration given by the order of the
operators from the filtration given by the power of $x$. It turns out that
on these algebras there are three `unique' trace functionals
$\Tr_{\partial,\sigma}$ on $\mathcal{A}_\sigma$,
$\mathcal{A}_\partial$ and $\mathcal{A}_{\partial,\sigma}$,
$\Tr_{\partial}$ on $\mathcal{I}_\partial$, and $\Tr_{\sigma}$ on
$\mathcal{I}_\sigma$. They are conic versions of the functionals
studied in \cite{Me-Ni}. $\Tr_{\partial,\sigma}$ and
$\Tr_{\sigma}$ were also considered in \cite{Sc97}. Their
definitions and basic properties are given in
Section~\ref{TraceFunctionals}.

In \cite{GiPHD98,GiHeat01}, the first author defines natural
conditions (called parameter-ellipticity) on a cone differential
operator $A \in x^{-m}\mathrm{Diff}^m_b(M)$ which ensure that the
heat operator $e^{-tA}$ exists. Moreover, a full asymptotic
expansion of $\Tr e^{-tA}$ was obtained. Later in
\cite{LoHeatc01}, this expansion was generalized to $\Tr Pe^{-tA}$
where $P$ is a cone differential operator. In the present work,
the techniques of \cite{LoHeatc01} are expanded to analyze the
generalized heat trace expansion when $P$ is pseudodifferential.

Let $P\in x^{-p} \Psi^{m'}_b(M)$, $p,m'\in\mathbb R$, be a
pseudodifferential cone operator. Assume that $p<m$ so that
$Pe^{-tA}$ is of trace class (on appropriate weighted Sobolev
spaces) for all $t>0$. In Section~\ref{HeatTrace}, we obtain an
asymptotic expansion of the form
\begin{equation*}
 \Tr Pe^{-tA} \sim_{t\to 0^+}\;
  \sum_{k=0}^\infty c_k\, t^{\xi_k} +
  \sum_{k=0}^\infty c_k' (\log t)\, t^{\eta_k} +
  \sum_{k=0}^\infty c_k'' (\log t)^2\, t^{\omega_k}
\end{equation*}
and give explicit formulas for all the $\log t$ coefficients
$c_k'$ and $c_k''$ in terms of the trace functionals mentioned
above, see Theorems~\ref{thm:main1}, \ref{thm:main2} and their
corollaries. In particular, the coefficient of $\log t$ is
\[ -\frac 1{m} \Tr_\sigma(P) -\frac 1{m}\Tr_\partial(P)
   -\frac 1{m^2} \Tr_{\partial,\sigma}(P), \]
and $\Tr_{\partial,\sigma}(P)= -m^2\;\times$ the coefficient of
$(\log t)^2$.

Our results rely on the parametrix construction within a suitable
parameter-dependent calculus which is presented in
Section~\ref{ParameterCalculus}, and on basic properties of the Laplace
and Mellin transforms discussed in Section~\ref{LMTransforms}.
Let us finally mention that the study of heat trace asymptotics on
manifolds with conic singularities was initiated by J.~Cheeger
\cite{Ch79,Ch83} and further developed by many other authors, a list
of references can be found in \cite{GiHeat01}, \cite{LoHeatc01}.

%%%%%%%%%%%%%%%%%%%%%%%%%%%%%%%%%%%%%%%%%%%%%%%%%%%%%%%%%%%%
\section{Operator algebras and trace functionals}\label{TraceFunctionals}
%%%%%%%%%%%%%%%%%%%%%%%%%%%%%%%%%%%%%%%%%%%%%%%%%%%%%%%%%%%%
Recall that an $n$-dimensional manifold with corners $M$ is a
manifold with atlas given by local models of the form
$[0,\infty)^k \times \R^{n-k}$, where $k$ can run anywhere between
$0$ and $n$. (This definition is suitable for our purposes, but is
a bit more general than the standard definition, cf.~\cite{Mel93}.)
For any $\alpha \in \R$, the $b$-alpha density
bundle $\Omega_b^\alpha$ is the line bundle on $M$ with local
basis of the form $|(dx/x)dy|^\alpha$ on a patch $[0,\infty)_x
\times \R^{n-k}_y$. We will either use $\alpha = 1/2$ or $\alpha =
1$; in the latter case, we will write $\Omega_b$ for $\Omega_b^1$.

\begin{figure} \centering
\setlength{\unitlength}{.30mm}
\begin{picture}(120,80)(-30,-10)
\thicklines \put(-21,20){\line(0,1){55}}
\put(10,-10){\line(1,0){55}} \qbezier(-21,20)(6,15)(9,-10)
\put(55,-8){$\mathit{rb}$} \put(-6,17){$\mathit{ff}$}
\put(-19,60){$\mathit{lb}$} {\thinlines \put(1,9){\line(1,1){63}}
} \put(60,60){$\Delta_b$}
%\put(6,2){\vector(2,1){30}}
%\put(36,15){$\rho$}
\put(-50,60){$M^2_b$}
\end{picture}
\caption{The manifold $M^2_b$ is $M^2$ ``blown-up'' at the corner
$Y^2$. The submanifold $\Delta_b$ is the the diagonal in $M^2$
lifted to $M^2_b$; that is, written in the polar coordinates of
$M^2_b$. %The function $\rho$ is a boundary defining
%function for $\mathit{ff}$.
} \label{fig:M2b}
\end{figure}
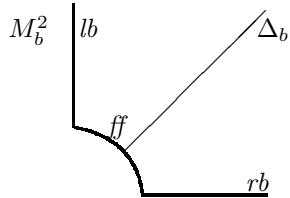

We now review the definition of $b$-pseudodifferential operators.
The standard reference is Melrose's book \cite{Mel93}. Let $M$ be
an $n$-dimensional compact manifold with connected boundary $Y =
\partial M$. Recall that $M^2_b$ is the manifold with corners that
has an atlas consisting of the usual coordinate patches on $M^2
\setminus Y^2$ together with polar coordinate patches over $Y^2$
in $M^2$. The typical picture of $M^2_b$ is shown in
Figure~\ref{fig:M2b}. The boundary hypersurfaces \emph{lb},
\emph{rb}, and \emph{ff} stand for ``left boundary'',
``right boundary'', and ``front face''.

For symmetry reasons, it is common practice to have
$b$-pseudodifferential operators act on $b$-half densities rather
than on functions. Given $m \in \R$, the space
$\Psi^m_b(M,\Omega_b^\frac12)$ consists of operators $A$ on
$C^\infty(M,\Omega_b^\frac12)$ that have a Schwartz kernel $K_A$
satisfying the following two conditions:
\begin{enumerate}
\item given $\varphi \in C^\infty(M^2_b \setminus \Delta_b)$,
we can write $\varphi K_A$ in the form $k \cdot \nu$, where $k$ is
a smooth function on $M^2_b$ vanishing to infinite order at
\emph{lb} and \emph{rb}, and where $\nu \in
C^\infty(M^2_b,\Omega^\frac12_b)$,
\item given a coordinate patch of $M^2_b$ overlapping $\Delta_b$
of the form $\mathcal{U}_y \times \R^n_\eta$ such that $\Delta_b
\cong \mathcal{U} \times \{0\}$, and given $\varphi \in
C^\infty_c(\mathcal{U} \times \R^n)$, we have \begin{equation}
\varphi K_A = \int e^{i\xi \cdot \eta}a(y,\xi) \dbar \xi \cdot
\nu,\label{KA} \end{equation} where $\nu \in
C^\infty(M^2_b,\Omega_b^\frac12)$, and where $a(y,\xi)$ is a
classical symbol of order $m$.
\end{enumerate}

Henceforth, we will fix a boundary defining function $x$ for $Y$.
Let $p$, $m \in \R$. Given $A \in x^{-p}\Psi^m_b(M,\Omega_b^\frac12)$, we
define its (Wodzicki) residue density $\res A$
% \in x^{-p}C^\infty(M,\Omega_b)$,
as follows, cf.~\cite{WoI87}. Write
the kernel of $A$ locally as in \eqref{KA} above. Observe that
since $\Delta_b \cong M$, the coordinate patch $\mathcal{U}$ can
be considered a coordinate patch on $M$. We define $\res A$
locally on the patch $\mathcal{U}$ by
\begin{equation}
\omega_y(A) := \int_{|\xi|=1} a_{-n}(y,\xi) \dbar \xi \cdot
\nu|_{\Delta_b},\label{resA}
\end{equation}
where $a_{-n}(y,\xi)$ is the homogeneous term of degree
$-n$ in the symbol expansion of $a(y,\xi)$. Note that if $m
\not\in \mathbb Z$, then $\res A = 0$ since $a(y,\xi)$ has no
homogeneous component of degree $-n$. Also note that
$\nu|_{\Delta_b} \in  C^\infty(M,\Omega_b)$ as one can check using
local coordinates. The local expression \eqref{resA} turns out to
be independent of coordinates; a nice proof of this fact can be
found in \cite{Fo-Go-Le-Sc96}. Thus, the local expressions
\eqref{resA} actually define a global density $\res A \in
x^{-p}C^\infty(M,\Omega_b)$.

We now review various trace functional introduced by Melrose and
Nistor \cite[Sec.~5]{Me-Ni}, and in a slightly different setting
by Schrohe \cite{Sc97}. We remark that Melrose and Nistor studied
cusp operators, although by continuity many of their results apply
to $b$-operators. In what follows, we will freely use their
results presented from the perspective we need them.
%%%%%%%%%%%%%%%%%%%%%%%%%%%%%%%%%%%%%%%%%%%%%%%%%%%%%%%%%%%%
\subsection*{The functional
$\mathbf{Tr}_{\boldsymbol{\partial},\boldsymbol{\sigma}}$}

Given $P\in x^{-p}\Psi^m_b(M,\Omega_b^\frac12)$, $p,m\in\R$,
we define
\begin{equation}
\Tr_{\partial,\sigma}(P):= \frac{1}{p!} \int_Y
\partial_x^{p} \{x^{p}\res P\}|_{x=0}. \label{Trps}
\end{equation}
If $p \not\in \N_0$, then $\Tr_{\partial,\sigma}(P)$ is defined to
be zero. The expression $\partial_x^{p} \{x^{p}\res P\}|_{x=0}$ is
defined as follows. Since $\res P \in x^{-p}
C^\infty(M,\Omega_b)$, in a collar $X \cong [0,1)_x \times Y$ near
$Y$, we can write $x^p \res P = f(x) \frac{dx}{x}$, where $f(x)$
is a smooth family of densities on $Y$. We define $\partial_x^{p}
\{x^{p}\res P\}|_{x=0}:= \partial_x^{p} f(0)$. It is not obvious
that $\Tr_{\partial,\sigma}(P)$ is defined independent of the
boundary defining function $x$. We will prove that it is
independently defined following Lemma \ref{lem:bint} below.

A more transparent way to view $\Tr_{\partial,\sigma}(P)$ is as
follows. Assuming that $p \in \N_0$, expand $\res P$ near the
boundary $Y$ in Taylor series: $\res P \sim \sum_{j=-p}^\infty
x^{j} \res{P_j}$. Then as $\frac{1}{p!}
\partial_x^{p} \{ x^{p}\res P\}|_{x=0} = \res{P_0}$, we have
\begin{equation}
\Tr_{\partial,\sigma}(P)= \int_Y \res{P_0}.\label{Trpsres}
\end{equation}
Thus, $\Tr_{\partial,\sigma}$ is a natural way to restrict the
Wodzicki residue trace to the boundary.

\begin{lemma} \label{lem:bint}
Let $u \in x^{-p}C^\infty(M,\Omega_b)$, $p \in\mathbb{R}$. Then
the function $\C \ni z \mapsto \int_M x^z u$ is well-defined for
$\Re z > p$ and it extends to be a meromorphic function on $\C$
with only simple poles at $z = p,\, p-1,\, p-2,\ldots$. In
particular, it has a pole at $z=0$ iff $p \in \N_0$ in which
case, its residue is given by $\frac{1}{p!} \int_Y
\partial_x^p \{x^p u\}|_{x=0}$. The regular value of the function
$\C \ni z \mapsto \int_M x^z u$ is
called the $b$-integral of $u$ and is denoted by $\trint{b}_M u$.
\end{lemma}
\begin{proof}
In a neighborhood $M \cong [0,1)_x \times Y$ of $Y$, we can write
$u = x^{-p}v(x)\frac{dx}{x}$ where $v(x)$ is a smooth family of
densities on $Y$. Expanding $v(x)$ in Taylor series at $x=0$ shows
that $u \sim \sum_{j=-p}^\infty x^j u_j$, where
\[ u_j = \frac{1}{(p+j)!} \partial_x^{p+j} v(0)
=\frac{1}{(p+j)!}\partial_x^{p+j} \{ x^p u\}|_{x=0},\;\; j= -p,
-p+1, \ldots .\] Since $\int_0^1 x^{z+j}\frac{dx}{x}=
\frac{1}{z+j}$ our lemma follows.
\end{proof}

\begin{lemma} \cite[Lem.~6]{Me-Ni} \label{lem:TrpsxA}
The functional $\Tr_{\partial,\sigma}$ is defined independent of
the choice of boundary defining function $x$.
\end{lemma}
\begin{proof}
Suppose that $x'$ is another boundary defining function. Then $x'
= ax$ where $0<a \in C^\infty(M)$. Denote by
$\Tr_{\partial,\sigma,x}$ and $\Tr_{\partial,\sigma,x'}$ the
functional $\Tr_{\partial,\sigma}$ defined using the boundary
defining functions $x$ and $x'$, respectively. Then by Lemma
\ref{lem:bint},
\[ \Tr_{\partial,\sigma,x}(P)-\Tr_{\partial,\sigma,x'}(P)
  =\Res_{1} \int_M  x^z P(z), \] where $\Res_{1}$ denotes the
residue at $z=0$, and where $P(z) = (1 -a^z)\res P$. Since $P(0) =
0$, by Lemma \ref{lem:bint} it follows that $\int_M x^z P(z)$ is
regular at $z=0$. Thus,
$\Tr_{\partial,\sigma,x}(P)=\Tr_{\partial,\sigma,x'}(P)$.
\end{proof}

We now consider in what sense the functional $\Tr_{\partial,\sigma}$
defines a trace. To do so, we need to introduce some operator
algebras. The main algebra is
\[ x^{-\mathbb{Z}}\Psi^{\mathbb{Z}}_b(M,\Omega_b^\frac12) =
  \bigcup_{p\in \mathbb{Z}} \bigcup_{m \in \mathbb{Z}} x^{-p}
  \Psi^{m}_b(M,\Omega_b^\frac12).\]
Consider now the following quotient algebras:
\begin{equation*}
\begin{split}
 \mathcal{A}_\sigma &:=
 x^{-\mathbb{Z}}\Psi^{\mathbb{Z}}_b(M,\Omega_b^\frac12)/
 x^{-\mathbb{Z}}\Psi^{-\infty}_b(M,\Omega_b^\frac12),\\
 \mathcal{A}_\partial &:=
 x^{-\mathbb{Z}}\Psi^{\mathbb{Z}}_b(M,\Omega_b^\frac12)/
 x^{\infty}\Psi^{\mathbb{Z}}_b(M,\Omega_b^\frac12),
\end{split}
\end{equation*}
\begin{equation*}
 \mathcal{A}_{\partial,\sigma} :=
 x^{-\mathbb{Z}}\Psi^{\mathbb{Z}}_b(M,\Omega_b^\frac12)/
 \{ x^{-\mathbb{Z}}\Psi^{-\infty}_b(M,\Omega_b^\frac12) +
    x^{\infty}\Psi^{\mathbb{Z}}_b(M,\Omega_b^\frac12)\}.
\end{equation*}
These algebras were first considered by Melrose and Nistor in
\cite{Me-Ni}, where the following theorem is also proved.

\begin{theorem}
The functional $\Tr_{\partial,\sigma}$ defines a trace on the
algebras $\mathcal{A}_\sigma$, $\mathcal{A}_\partial$, and
$\mathcal{A}_{\partial,\sigma}$, and is the unique such trace in
the sense that any other trace on any of these algebras is a
constant multiple of $\Tr_{\partial,\sigma}$.
\end{theorem}

%%%%%%%%%%%%%%%%%%%%%%%%%%%%%%%%%%%%%%%%%%%%%%%%%%%%%%%%%%%%
\subsection*{The functional
$\mathbf{Tr}_{\boldsymbol{\partial}}$}

Fix a holomorphic family $Q(z) \in
x^{\alpha z}\Psi^{\beta z}_b(M,\Omega_b^\frac12)$,
$\alpha,\beta\in\R$, of operators such that $Q(0)=\mathrm{Id}$.
Let $P \in x^{-p}\Psi^{m}_b(M,\Omega_b^\frac12)$, $p,m\in \R$, and
denote by $(P Q(z))|_\Delta$ the restriction of the Schwartz kernel
of $PQ(z)$ to the diagonal $\Delta$ in $M^2$. Then by
\cite[Lem.~4]{Me-Ni}, $(PQ(z))|_\Delta$ defines a meromorphic function
on $\C$, taking values in $x^{\alpha z-p}C^\infty(M,\Omega_b)$
(using that $\Delta \cong M$), with possible simple poles at those
$z\in\C$ with $\beta z =-n-m+k$ for $k= 0,1, \ldots$. In particular,
\[ \Res_{0} (P Q(z))|_\Delta \in x^{-p} C^\infty(M,\Omega_b) \]
is well defined, where $\Res_{0}$ means ``the regular value at $z=0$''.
We define
\begin{equation}
\Tr_{\partial}(P):= \frac{1}{p!} \int_Y \partial_x^{p} \{ x^{p}
\Res_{0} (P Q(z))|_\Delta\}|_{x=0}.\label{TrpP}
\end{equation}
If $p \not\in \N_0$, then $\Tr_{\partial}(P)$ is defined to be
zero. The functional $\Tr_\partial$ was first introduced in
\cite{Me-Ni}. The same argument found in Lemma~\ref{lem:TrpsxA}
shows that $\Tr_{\partial}$ is defined independent of the boundary
defining function chosen. (Unfortunately, $\Tr_\partial$ does
depend on the regularizing operator $Q(z)$; however, its
dependence on $Q(z)$ can be explicitly determined, see
\cite[Lem.~11]{Me-Ni}.)

For $P$ of sufficiently large negative order, $\Tr_\partial(P)$
has a natural interpretation. Indeed, if $m<-n$, then observe that
$(P Q(z))|_\Delta$ is regular at $z=0$ with value $K_P|_\Delta \in
x^{-p}C^\infty(M,\Omega_b)$. Assume that $p \in \N_0$ and expand
$K_P|_\Delta$ in Taylor series at $x = 0$: $K_P|_\Delta \sim
\sum_{j=-p}^\infty x^j K_{P,j}|_\Delta$. Then as $\partial_x^{p}
\{ x^{p} K_{P}|_\Delta\}|_{x=0} = K_{P,0}|_\Delta$, we have
\[ \Tr_\partial (P) = \int_Y K_{P,0}|_\Delta.\]
Thus, for $P$ of sufficiently negative order, $\Tr_\partial (P)$
is a type of $L^2$-trace of $P$ restricted to the boundary of $M$.

Consider now the algebra \[ \mathcal{I}_\partial  :=
 x^{-\mathbb{Z}}\Psi^{-\infty}_b(M,\Omega_b^\frac12)/x^\infty \Psi^{-\infty}_b
 (M,\Omega_b^\frac12).\]  The following theorem is proved in
\cite{Me-Ni}.
\begin{theorem}
The functional $\Tr_{\partial}$ defines a trace on
$\mathcal{I}_\partial$ and is the unique such trace in the sense that
any other trace on $\mathcal{I}_\partial$ is a constant multiple
of $\mathrm{Tr}_{\partial}$.
\end{theorem}

%%%%%%%%%%%%%%%%%%%%%%%%%%%%%%%%%%%%%%%%%%%%%%%%%%%%%%%%%%%%
\subsection*{The trace functional
$\mathbf{Tr}_{\boldsymbol{\sigma}}$}

Given $P \in x^{-p}\Psi^{m}_b(M,\Omega_b^\frac12)$, $p,m \in \R$,
we define
\begin{equation}
\Tr_{\sigma}(P):= \rint{b}_M \res P, \label{resP}
\end{equation}
where $\trint{b}$ is the $b$-integral introduced in
Lemma~\ref{lem:bint}. This functional was first introduced in
\cite{Sc97} when $P$ vanishes at the boundary (in which case,
$\trint{b}_M \res P = \int_M \res P$ is just the usual integral of
$\res P$ over $M$). In the generality presented in \eqref{resP},
$\Tr_\sigma$ was first studied in \cite{Me-Ni}. $\Tr_\sigma$ is
the natural generalization of the Wodzicki residue trace to cone
pseudodifferential operators.

To see how the functional $\Tr_{\sigma}$ depends on the boundary
defining function $x$, we first prove the following lemma.

\begin{lemma} \label{lem:bintxx'}
Let $x' = ax$ be another boundary defining function for $Y$ where
$0<a \in C^\infty(M)$ and let $u\in x^{-p}C^\infty(M,\Omega_b)$,
$p \in \R$. Denote by $\trint{b,x}_M u$ and $\trint{b,x'}_M u$ the
$b$-integral of $u$ as defined by using the boundary defining
function $x$ and $x'$, respectively. Then,
\begin{equation}
\rint{b,x}_M u = \rint{b,x'}_M u + \frac{1}{p!} \int_Y
\partial_x^{p}\{x^p \log a\, u\}|_{x=0}.\label{bintxx'formula}
\end{equation}
If $p \not\in \N_0$, then the last term is understood to be equal
to zero.
\end{lemma}
\begin{proof}
We can write $(x')^z u - x^z u = z x^z a(z)u$, where $a(z) = (a^z - 1)/z$.
As $a(0) = \log a$, we have
\[ \rint{b,x'}_M u\; - \;\rint{b,x}_M u = \Res_0 z\hspace{-.3em} \int_M x^z a(z)
u = \Res_{1}\int_M x^z \log a\, u,\] where $\Res_{1}$ denotes the
residue at $z=0$. \eqref{bintxx'formula} now follows from
Lemma~\ref{lem:bint}.
\end{proof}

\begin{lemma}\emph{\cite[Lem.~9]{Me-Ni}}
Given another boundary defining function $x' = ax$ where $0 < a
\in C^\infty(M)$, the trace functionals $\Tr_\sigma$ defined with
respect to $x'$ and $x$, respectively, are related by
\[ \Tr_{\sigma,x'}(P) - \Tr_{\sigma,x}(P)
  =\Tr_{\partial,\sigma}(\log a\, P). \]
\end{lemma}
\begin{proof}
By Lemma \ref{lem:bintxx'}, we have
\begin{align*}
\Tr_{\sigma,x'}(P) - \Tr_{\sigma,x}(P)
 &=\rint{b,x}_M \res P  -\rint{b,x'}_M \res P\\
 &=\frac{1}{p!} \int_Y \partial_x^{p}
   \{ x^{p}\res{\log a\, P}\}|_{x=0}\\
 &=\Tr_{\partial,\sigma}(\log a\, P).
\end{align*}
\end{proof}
The final algebra we consider is
\[ \mathcal{I}_\sigma :=
 x^{\infty}\Psi^{\mathbb{Z}}_b(M,\Omega_b^\frac12)/x^\infty \Psi^{-\infty}_b
 (M,\Omega_b^\frac12). \]
The following theorem is proved in \cite{Me-Ni}.

\begin{theorem}
The functional $\Tr_{\sigma}$ defines a trace on
$\mathcal{I}_\sigma$ and is the unique such trace in the sense
that any other trace on $\mathcal{I}_\sigma$ is a constant
multiple of $\Tr_{\sigma}$.
\end{theorem}

%%%%%%%%%%%%%%%%%%%%%%%%%%%%%%%%%%%%%%%%%%%%%%%%%%%%%%%%%%%%
\section{Parameter-dependent calculus}\label{ParameterCalculus}
%%%%%%%%%%%%%%%%%%%%%%%%%%%%%%%%%%%%%%%%%%%%%%%%%%%%%%%%%%%%

The space $x^{-m}\Diff_b^m(M,\Omega_b^{\frac 1{2}})$, $m>0$, is the space
of differential operators  on $M$ that near the boundary, take the form
\[ A=x^{-m} \hspace{-.8em}
\sum_{k+|\alpha|\le m} a_{k\alpha}(x,y)(xD_x)^k D_y^\alpha \]
with coefficients smooth up to $x=0$. To such an $A$ we associate the
operator
\[ A_0=x^{-m} \hspace{-.8em}
\sum_{k+|\alpha|\le m} a_{k\alpha}(0,y)(xD_x)^k D_y^\alpha \]
by freezing the coefficients at $x=0$.

\begin{definition}
Let $A\in x^{-m}\Diff_b^m(M,\Omega_b^{\frac 1{2}})$ and let
$\Lambda$ be a sector in $\C$ containing the origin. The operator
family $A-\lambda$ is said to be {\it parameter-elliptic} with
respect to $\alpha\in\R$ on $\Lambda$ if and only if
\begin{enumerate}
\item $\sigma^m_{\psi,b}(A)(\xi)-\lambda$ is invertible for all $\xi\not=0$
and $\lambda\in\Lambda$,
\item $\mathrm{spec}_b(A)\cap \{z\in\C\,|\; \Im z=-\alpha\}=\varnothing$,
\item $A_0-\lambda: x^{\alpha}H_{b}^s(Y^\wedge,\Omega_b^{\frac 1{2}})
\to x^{\alpha-m}H_{b}^{s-m}(Y^\wedge,\Omega_b^{\frac 1{2}})$ is invertible
for every $\lambda\in\Lambda$ sufficiently large, and for some $s\in\R$.
\end{enumerate}
Here, cf.\ \cite[Sec.~3]{GiHeat01}, $\sigma^m_{\psi,b}(A)(\xi)$ is the
totally characteristic principal symbol of $x^m A$,
$\mathrm{spec}_b(A)$ denotes the boundary spectrum of $x^m A$,
$Y^\wedge:=Y\times\overline\R_+$, and $H_b^s$ denotes the
totally characteristic Sobolev space of order $s$.
\end{definition}

The following proposition is proved in \cite[Th.~3.2]{GiHeat01}.
\begin{proposition}
If $A-\lambda$ is parameter-elliptic with respect to $\alpha$ on
$\Lambda$, then
\[ A-\lambda: x^{\alpha}H_{b}^s(M,\Omega_b^{\frac 1{2}})
  \to x^{\alpha-m}H_{b}^{s-m}(M,\Omega_b^{\frac 1{2}}) \]
is invertible for every $\lambda\in\Lambda$ sufficiently large,
and all $s\in\R$.
\end{proposition}

We next analyze the resolvent $(A -\lambda)^{-1}$ within a
suitable parameter-dependent pseudodifferential calculus. We begin
by describing the parameter-dependent symbols as defined in
\cite{LoHeatc01}. Related symbol classes can be found in
\cite{GrSe95}.

For $m$, $p\in\R$ and $d>0$ we define $S^{m,p,d}(\R^n;\Lambda)$ as
the space of functions $a\in C^\infty(\R^n\times \Lambda)$ such that
\[ |\partial_\xi^\alpha \partial_\lambda^\beta a(\xi,\lambda)|
 \le C_{\alpha \beta} (1+|\xi|)^{m-p-|\alpha|}
 (1+|\xi|+|\lambda|^{1/d})^{p-d|\beta|}. \]
The space $S_r^{m,p,d}(\R^n;\Lambda)$, $p/d\in\mathbb Z$, consists
of elements $a\in S^{m,p,d}(\R^n;\Lambda)$ such that if we set
\[ \tilde a(\xi,z):= z^{p/d}a(\xi,1/z), \]
then $\tilde{a} (\xi,z)$ is smooth at $z=0$, and
\[ |\partial_\xi^\alpha \partial_z^\beta \tilde a(\xi,z)|
 \le C_{\alpha \beta} (1+|\xi|)^{m-p-|\alpha|+d|\beta|}
 (1+|z||\xi|^{d})^{p/d-|\beta|} \]
uniformly for $|z|\le 1$. Let $S_{r,c\ell}^{m,p,d}(\R^n;\Lambda)$
be the class of elements $a\in S_r^{m,p,d}(\R^n;\Lambda)$ such that for
every $N\in\N$,
\begin{equation} a(\xi,\lambda)=\sum_{j=0}^{N-1}
\chi(\xi)a_{m-j}(\xi,\lambda)
  + r_{N}(\xi,\lambda) \label{symbexp} \end{equation}
with $r_N\in S_r^{m-N,p,d}(\R^n;\Lambda)$, $a_{m-j}$ anisotropic
homogeneous of degree $m-j$, i.e.,
\[ a_{m-j}(\delta\xi,\delta^d\lambda)=\delta^{m-j}a_{m-j}(\xi,\lambda)
  \;\text{ for } \delta>0, \]
and $z^{p/d} a_{m-j}(\xi,1/z)$ smooth at $z=0$. Finally, an
element $a\in S_{r,c\ell}^{m,p,d}(\R^n;\Lambda)$ is said to be
{\em holomorphically tempered} if it is holomorphic in a
neighborhood of $\Lambda$, and there exists an $\eps>0$ such that
each $a_{m-j}(\xi,\lambda)$ is holomorphic in the region
\[ \big\{(\xi,\lambda)\in (\R^n \setminus \{0\})\times \C : \lambda \in
\Lambda\ \text{or\ } |\lambda|  \le \eps |\xi|^d\ \text{or\
}\tfrac{1}{\eps} |\xi|^d \le |\lambda|\big\}. \]

We now define our corresponding spaces of parameter-dependent cone
operators, cf.\ \cite{LoRes01,LoRes201}. Let $\rho$ denote a
boundary defining function for the front face \emph{ff} of
$M^2_b$. Our first space
$\Psi_{c,\Lambda}^{m,p,d}(M,\Omega_b^{\frac 1{2}})$ is defined as
the space of those operator families $A(\lambda)\in
C^\infty(\Lambda,\Psi_b^m(M,\Omega_b^{\frac 1{2}}))$ whose
Schwartz kernels $K_{A(\lambda)}$ satisfy the following two
conditions:
\begin{enumerate}
\item given $\varphi\in C_c^\infty(M_b^2\setminus \Delta_b)$, then
$\varphi K_{A(\lambda)}=R(\rho^d\lambda)$, where $R(\lambda)$ is
smooth in $\lambda$ taking values in
$\Psi^{-\infty}_b(M,\Omega_b^\frac12)$ and satisfies the estimates
\[ |\partial_\lambda^\alpha R(\lambda)|\leq
  C_\alpha (1 + |\lambda|)^{p/d - |\alpha|},\ \mbox{for all\ }\alpha
  \in \mathbb{N}_0, \;\lambda\in\Lambda, \]
\item given a coordinate patch on $M_b^2$ overlapping $\Delta_b$
of the form $\mathcal{U}_y \times
\mathbb{R}^n_\eta$ such that $\Delta_b \cong \mathcal{U}
\times\{0\}$, and given $\varphi \in C_c^\infty(\mathcal U\times
\R^n)$, then
 \[ \varphi K_{A(\lambda)} = \int e^{i\xi\cdot \eta}
 a(y,\xi,\rho^d\lambda)\dbar\xi \cdot \nu, \]
where $\nu \in C^\infty(M^2_b,\Omega_b^\frac12)$, and where
\[ y\mapsto a(y,\xi,\lambda)\in
  C^\infty(\mathcal{U},S_{r,c\ell}^{m,p,d}(\R^n;\Lambda)).\]
\end{enumerate}
Furthermore, $A(\lambda)$ is holomorphically tempered if it is
holomorphic in a neighborhood of $\Lambda$ and its corresponding
local symbols are holomorphically tempered.

We sketch the definition of our next parameter-dependent space;
see \cite{LoRes01} for the precise definition. Assume that a branch 
of $\log$ is defined and fixed on $\Lambda$, and let 
\[ \Lambda^{1/d}:= \{\lambda^{1/d}\,|\, \lambda\in\Lambda\setminus 0\}
   \cup \{0\}. \]
Let $\mathcal{T}$ denote the manifold $M^2_b \times \Lambda^{1/d}$ 
blown-up at the submanifold $\mathit{ff}\times 
\{\lambda^{1/d}=\infty\}$. The boundary hypersurface of $\mathcal{T}$
originating from $M^2_b \times \{\lambda^{1/d} = \infty\}$ is called 
the ``boundary at $\infty$''. We define
$\Psi^{-\infty,d,\mathcal{E}}_{c,\Lambda}(M,\Omega_b^\frac12)$ 
as the class of operators whose Schwartz
kernels can be written in the form $k \cdot \nu$, where $\nu \in
C^\infty(M^2_b,\Omega_b^\frac12)$, and where $k$ is a function on
$\mathcal{T}$ which vanishes to infinite order at the boundary at
$\infty$, and is of asymptotic type $\mathcal{E}$ at the rest of
the faces of $\mathcal{T}$. Here, $k$ is of asymptotic type
$\mathcal{E}$ means, roughly speaking, that if $H$ is a boundary
hypersurface of $\mathcal{T}$ (except the boundary at $\infty$),
then $\mathcal{E}$ associates to $H$ a set of numbers $E_H$ such
that $k$ can be expanded at $H$ in powers and powers of the
logarithm of the corresponding boundary defining function; the
powers that may occur in the expansion are given exactly by the
set $E_H$.

Our final space of operators is
$\Psi^{-\infty,\mathcal{F}}_\Lambda(M,\Omega_b^\frac12)$, where
$\mathcal{F}=(F_1,F_2)$ is an index family for $M^2$. An element
$S(\lambda)\in \Psi^{-\infty,\mathcal{F}}_\Lambda(M,\Omega_b^\frac12)$
is a family of operators whose Schwartz kernels $K_{S(\lambda)}$
can be written in the form $k\cdot \nu$, where $\nu \in
C^\infty(M^2,\Omega_b^\frac12)$, and
where $k$ is a function on $M^2 \times \Lambda$ that vanishes to
infinite order at $\lambda = \infty$, and is of asymptotic type
$\mathcal{F}$ on $M^2$.

The next theorem follows from \cite[Th.~2.16]{LoHeatc01} and from
the proof of \cite[Th.~3.1]{LoRes201}.
\begin{theorem}\label{GenResExpansion}
Let $P\in x^{-p}\Psi_b^{m'}(M,\Omega_b^{\frac 1{2}})$, $p$,
$m'\in\R$. Let $A\in x^{-m}\Diff_b^m(M,\Omega_b^{\frac 1{2}})$ be
such that $A-\lambda$ is parameter-elliptic with respect to
$\alpha\in\R$ on a sector $\Lambda$. Then we get a decomposition
\begin{equation*}
 P(A-\lambda)^{-1} = Q(\lambda) + R(\lambda) + S(\lambda),
\end{equation*}
where $Q(\lambda) \in
x^{m-p}\Psi_{c,\Lambda}^{m'-m,-m,m}(M,\Omega_b^{\frac 1{2}})$,
$R(\lambda)\in x^{-p}\Psi_{c,\Lambda}^{-\infty,m,
\mathcal G(\alpha)}(M,\Omega_b^{\frac1{2}})$, and $S(\lambda)\in
x^{-p}\Psi_{\Lambda}^{-\infty,\mathcal F(\alpha)}(M,\Omega_b^{\frac 1{2}})$
for some index families $\mathcal{G}(\alpha)$ and $\mathcal F(\alpha)$.
$Q(\lambda)$ being holomorphically tempered.
\end{theorem}
\begin{remark}
The index families $\mathcal G(\alpha)$ and $\mathcal F(\alpha)$
are given explicitly in terms of the boundary spectrum
$\mathrm{spec}_b(A)$. For their precise definition we refer to
\cite[Section~3.2]{LoRes01}. They will play no role in the trace
expansions that we are looking at, so we do not need their
explicit descriptions here.
\end{remark}

%%%%%%%%%%%%%%%%%%%%%%%%%%%%%%%%%%%%%%%%%%%%%%%%%%%%%%%%%%%%
\section{Analysis of the Laplace and Mellin transforms}\label{LMTransforms}
%%%%%%%%%%%%%%%%%%%%%%%%%%%%%%%%%%%%%%%%%%%%%%%%%%%%%%%%%%%%

In order to give a precise representation of the coefficients in
the asymptotic expansion of $\Tr P e^{-tA}$ we first want to
analyze the Laplace transform of a holomorphically tempered family
$Q(\lambda)\in x^{m-p}\Psi^{m'-m,-m,m}_{c,\Lambda}
(M;\Omega_b^\frac12)$, $m'$, $p\in\R$, $m>0$, $m >p$. Here, we
assume that $\Lambda$ is a sector of the form
\begin{equation}\label{sector}
\{ \lambda \in \C: \eps \le \arg \lambda \le 2 \pi - \eps\} \quad
\text{for some } 0< \eps < \pi/2.
\end{equation}

Since $Q(\lambda)$ is holomorphically tempered, there exists
$\delta>0$ such that $Q(\lambda)$ is holomorphic on and to the
left of the contour \[\Upsilon := \{ \lambda \in \C:\lambda\in
\partial \Lambda\ \text{for\ }|\lambda| \ge \delta,\ \text{or\ }
|\lambda|= \delta\ \text{for\ } 2\pi - \eps \le \arg \lambda \le
\eps\}.\] For $t>0$ define
\begin{equation}\label{IntegralQ}
 Q(t)= \frac i{2\pi} \int_\Upsilon e^{-t\lambda} Q(\lambda) d\lambda.
\end{equation} Then $Q(t)$ is an
operator of trace class on $x^{\alpha-m}H_b^s(M,\Omega_b^\frac12)$
for every $s\in\R$. This relies on the fact that for each $t>0$,
the kernel of $Q(t)$  is a smooth function on $M^2$ vanishing to
infinite order at $\partial M^2$, cf. \cite[Th.~3.21]{LoHeatc01}.

Let $M^2_b \cong [0,1)_x\times Y\times \R^n_\eta$ near the face
\emph{ff} with $\Delta_b \cong [0,1)_x\times Y\times \{0\}$, see
Figure \ref{fig:M2b}. For simplicity, in the analysis that follows
we assume that the kernel of $Q(\lambda)$ is supported on this
coordinate patch. Furthermore, since the variables on $Y$ enter in
the analysis more or less as parameters (cf.\ the proof of
\cite[Th.~3.21]{LoHeatc01}), we also omit the variables on $Y$.
Thus, we can write
\begin{equation} \label{KernelQ}
K_{Q(\lambda)}= \int e^{i\xi \cdot \eta} q(x,\xi,x^m\lambda)\dbar\xi
 \;\left|\tfrac{dx}{x}d\eta\right|^{1/2},
\end{equation}
where $q(x,\xi,\lambda)\in S_{r,c\ell}^{m'-m,-m,m}(\R^n;\Lambda)$
and $q(x,\xi,\lambda)=O(x^{m-p})$. Thus, \pagebreak
\begin{gather*}
 K_{Q(t)}= \int e^{i\xi\cdot \eta}\mathcal L_c(q)(t,x,\xi)
  \dbar\xi\;\left|\tfrac{dx}{x}d\eta\right|^{1/2}, 
\intertext{where}
 \mathcal L_c(q)(t,x,\xi)= \frac i{2\pi}
  \int_\Upsilon e^{-t\lambda} q(x,\xi,x^{m}\lambda) d\lambda. 
\end{gather*}
Finally,
\begin{equation}\label{TraceQLaplace}
 \Tr Q(t)= \iint \mathcal L_c(q)(t,x,\xi)\dbar\xi \,\tfrac{dx}{x}.
\end{equation}
Motivated by the relation between the heat trace and the zeta
function via the Mellin transform, we define
\begin{equation}\label{B(z)}
 B(z):= \frac 1{\Gamma(z)}\int_0^\infty t^{z-1}\Tr Q(t)dt.
\end{equation}
Thus, by means of the inverse Mellin transform we get
\begin{equation}\label{TraceQMellin}
 \Tr Q(t) = \frac 1{2\pi i}\int_{\Re z=\delta}t^{-z}B(z)\Gamma(z) dz
\end{equation}
for any sufficiently large $\delta\in\R$. The following
proposition (which is just an application of Cauchy's theorem)
suggests that the asymptotic expansion of $\Tr Q(t)$ as $t\to 0^+$
is determined by the poles of $B(z)\Gamma(z)$.

\begin{proposition}\label{InverseMellin}
Suppose that $\psi(z)$ is meromorphic with a single pole of order
$\ell +1$ at $z=w$. Let
\[ u(t)=\frac 1{2\pi i}\int_\gamma t^{-z} \psi(z)dz, \]
where $\gamma$ is a closed curve around $w$. Then, we can write
\[ u(t)=\sum_{k=0}^{\ell} \frac 1{k!} (-1)^{k}\, t^{-w}( \log t)^k r_{k+1}, \]
where $r_1, r_2,\dots$ are the coefficients of the Laurent
expansion of $\psi$.
\end{proposition}

Now, our goal is to determine all the poles of $B(z)\Gamma(z)$,
then push the contour $\Re z = \delta$ in \eqref{TraceQMellin}
down to $\Re z = -\infty$. Every time we pass a pole of
$B(z)\Gamma(z)$ we pick up a contour integral of $B(z)\Gamma(z)$
around that pole, which by Proposition~\ref{InverseMellin},
contributes powers and powers of the logarithm of $t$ to the
expansion of $\Tr Q(t)$.

We first determine the poles of $B(z)$. To do so, we want to write
it in terms of $q(x,\xi,\lambda)$. By \eqref{TraceQLaplace} and \eqref{B(z)},
and the fact that $\lambda^{-z}= \frac{1}{\Gamma(z)}
\int_0^\infty t^{z -1} e^{-t\lambda} dt$, we have
\begin{align*}
 B(z) &= \frac 1{\Gamma(z)}\int_0^\infty t^{z-1} \Tr Q(t)dt\\
   &= \frac 1{\Gamma(z)}\int_0^\infty t^{z-1}
      \iint \mathcal L_c(q)(t,x,\xi)\dbar\xi \,\tfrac{dx}{x} dt\\
   &= \iint \left(\frac 1{\Gamma(z)}\int_0^\infty t^{z-1}
      \mathcal L_c(q)(t,x,\xi)dt \right)\dbar\xi \,\tfrac{dx}{x}\\
   &= \iint \left( \frac i{2\pi} \int_\Upsilon \lambda^{-z}
      q(x,\xi,x^m\lambda) d\lambda\right)\dbar\xi \,\tfrac{dx}{x}\\
   &= \iint \hat q(x,\xi,z)\dbar\xi \,\tfrac{dx}{x},
\end{align*}
where
\begin{align*}
\hat q(x,\xi,z) &:= \frac i{2\pi} \int_\Upsilon
 \lambda^{-z} q(x,\xi,x^m\lambda) d\lambda\\
 &= x^{mz - m}\frac i{2\pi} \int_\Upsilon \lambda^{-z} q(x,\xi,\lambda)
 d\lambda\, \quad (\lambda \to x^{-m} \lambda).
  \end{align*}
According to \eqref{KernelQ}, $\hat q (x,\xi,z)$ is a local symbol
of the operator
\[ \M(Q)(z):= \frac i{2\pi} \int_\Upsilon \lambda^{-z}
  Q(\lambda) d\lambda.\]
By the definition of the space $S^{m'-m,-m,m}_{r, c\ell}(\R^n; \Lambda)$,
cf.~\eqref{symbexp}, a straight forward computation shows that
$\hat q(x,\xi,z)$ is a classical symbol of order $-mz + m'$. Now let us write
\begin{equation}\label{MellinQ}
  B(z)=\int_M \M(Q)(z)|_\Delta \quad\text{with}\quad
 \M(Q)(z)|_\Delta = \int_{\R^n} \hat q(x,\xi,z)\dbar\xi\; \tfrac{dx}{x}.
\end{equation}

As in Theorem 3.21 of  \cite{LoHeatc01}, $\M(Q)(z)|_\Delta$
extends to be a meromorphic function on $\C$ having simple poles
at $z=-z_k$, $z_k=\frac{k-m'-n}{m}$, $k\in\N_0$, with residues
\begin{equation*}
 \Res_1 \M(Q)(-z_k)|_\Delta = \frac 1{m}\res{\M(Q)(-z_k)},
\end{equation*}
where $\res{\M(Q)(-z_k)}$ is the corresponding Wodzicki residue
density, i.e.,
\[ \res{\M(Q)(-z_k)} = \int_{|\xi|=1}
   \hat q_{-n}(x, \xi,-z_k)\dbar\xi\; \tfrac{dx}{x}, \]
where $\hat q_{-n}(x, \xi,-z_k)$ denotes the homogeneous component
of order $-n$ in the expansion of $\hat q(x, \xi,-z_k)$. Note that
$\hat q(x,\xi,-z_k)$ is a classical symbol of order $k-n\in \mathbb Z$
since $\hat q (x,\xi,z)$ is classical of order $-mz + m'$.

Now it follows from Lemma \ref{lem:bint} that if $u\in x^{-p}
C_c^\infty(\overline\R_+)$, then $\int_0^\infty x^{mz} u(x)
\frac{dx}{x}$ has poles at $z=-\frac{j-p}{m}$, $j \in \N_0$ with
corresponding residues $\frac 1{m}\frac 1{j!}\partial_x^j(x^p
u)|_{x=0}$. Thus, as $\M(Q)(z)|_\Delta \in
x^{mz-p}C^\infty_c(\overline\R_+,\Omega_b)$, $B(z)$ has the poles
of $\M(Q)(z)|_\Delta$ and additional poles at $z=-\frac{j-p}{m}$.
On the other hand, $\Gamma(z)$ has simple poles at $z=-\ell$,
$\ell\in\N_0$, with residues given by $\frac{(-1)^\ell}{\ell!}$.
Altogether we get that $B(z)\Gamma(z)$ may have poles at
\begin{alignat}{2} \notag
 z &= -\ell &\text{for } &\ell=0,1,2,\dots \\ \label{BGPoles}
 z &= -\frac{j-p}{m} &\text{for } &j=0,1,2,\dots \\ \notag
 z &= -\frac{k-m'-n}{m} \quad &\text{for } &k=0,1,2,\dots
\end{alignat}
These are simple, double or triple order poles depending on $j$,
$k$ and $\ell$. In any case, only double and triple order poles
may produce $\log$ terms in the asymptotic expansion of $\Tr
Q(t)$, cf. Proposition~\ref{InverseMellin}. In summary, the poles
of $B(z)\Gamma(z)$ arise because of the poles of $\Gamma(z)$,
integrating $\M(Q)(z)|_\Delta$ in $x$, and from the poles of
$\M(Q)(z)|_\Delta$ itself. Using this information, we now write
down all the possible combinations for the higher order poles of
$B(z)\Gamma(z)$.

%%%%%%%%%%%%%%%%%%%%%%%%%%%%%%%%%%%%%%%%%%%%%%%%%%%%%%%%%%%%
\subsection*{Second order poles of $\mathbf{B(z)\Gamma(z)}$}
%%%%%%%%%%%%%%%%%%%%%%%%%%%%%%%%%%%%%%%%%%%%%%%%%%%%%%%%%%%%
Suppose that this function have a pole at $z=\zeta$.
According to \eqref{BGPoles} there are three cases where $\zeta$
may be a double pole.\\[1ex]
\underline{\it Case 1:}\quad
$\zeta=-\frac{j-p}{m}=-\frac{k-m'-n}{m}\not\in -\N_0$.
Then $j=p-\zeta m$, and the second order residue at $\zeta=-z_k$ is
\begin{align*}
 r_2 &= \Gamma(\zeta)
 \frac 1{mj!} \partial_x^j\Big\{x^{p - \zeta m }\frac 1{m}
 \res{\M(Q)(\zeta)} \Big\}\Big|_{x=0}\\
 &= \frac{\Gamma(-z_k)}{m^2(p+z_k m)!}
 \partial_x^{p+z_k m} \Big\{ x^{p +z_k m }
  \res{\M(Q)(-z_k)} \Big\}\Big|_{x=0}.
\end{align*}
\underline{\it Case 2:}\quad
$\zeta=-\ell=-\frac{k-m'-n}{m}\in -\N_0\;$ and $\zeta \not= - \frac{j-p}{m}$
for any $j$. Then
\begin{equation*}
% r_2 = \frac{(-1)^\ell}{\ell!}\;
% \rint{b} \frac 1{m}\res_x \M(Q)(-\ell)
 r_2 = \frac{(-1)^\ell}{m\ell!}\;\rint{b} \res{\M(Q)(-\ell)}.
\end{equation*}
\underline{\it Case 3:}\quad
$\zeta=-\ell=-\frac{j-p}{m}\;$ and $\zeta\not=-\frac{k-m'-n}{m}$ for any $k$.
In this case,
\begin{align*}
 r_2 &= \frac{(-1)^\ell}{\ell!}
 \frac 1{mj!}\partial_x^j\Big\{x^{j}\,
 \Res_0 \M(Q)(-\ell)|_{\Delta} \Big\}\Big|_{x=0}\\
 &= \frac{(-1)^\ell}{m\ell!(p+\ell m)!}
 \partial_x^{p+\ell m}\Big\{x^{p + \ell m}\,
 \Res_0 \M(Q)(-\ell)|_{\Delta}\Big\}\Big|_{x=0},
\end{align*}
where $\Res_0 \M(Q)(-\ell)|_{\Delta}$ means the regular value of
the kernel of $\M(Q)(z)$ restricted to the diagonal $\Delta$ in
$M^2$, evaluated at $z=-\ell$.

%%%%%%%%%%%%%%%%%%%%%%%%%%%%%%%%%%%%%%%%%%%%%%%%%%%%%%%%%%%%
\subsection*{Third order poles of $\mathbf{B(z)\Gamma(z)}$}
%%%%%%%%%%%%%%%%%%%%%%%%%%%%%%%%%%%%%%%%%%%%%%%%%%%%%%%%%%%%
They may only occur when in \eqref{BGPoles}
\[ \ell=\frac{j-p}{m}=\frac{k-m'-n}{m}
  \;\text{ for some } j,k,\ell\in\N_0. \]
In this case, the third order residue is given by
\begin{align*}
 r_3 &= \frac{(-1)^\ell}{\ell!}
 \frac 1{mj!} \partial_x^j\Big\{x^{j}\frac 1{m}
 \res{\M(Q)(-\ell)} \Big\}\Big|_{x=0}\\
 &= \frac{(-1)^\ell}{m^2\ell! (p+\ell m)!}
 \partial_x^{p+\ell m}\Big\{x^{p + \ell m}\res{\M(Q)(-\ell)} \Big\}\Big|_{x=0}.
\end{align*}
The second order residue can be written as a sum of the following
three expressions:
\begin{align*}
 r_{2,1} &= \frac{\Gamma_0(-\ell)}{m^2(p+\ell m)!}
 \partial_x^{p+\ell m}\Big\{ x^{p +\ell m }
  \res{\M(Q)(-\ell)} \Big\}\Big|_{x=0},\\
\intertext{where $\Gamma_0(-\ell)$ denotes the regular value of
$\Gamma(-\ell)$,}
 r_{2,2} &= \frac{(-1)^\ell}{m\ell!}\;
 \rint{b} \res{\M(Q)(-\ell)},\\
 r_{2,3} &= \frac{(-1)^\ell}{m\ell!(p+\ell m)!}
 \partial_x^{p+\ell m}\Big\{x^{p + \ell m}\,
 \Res_0 \M(Q)(-\ell)|_{\Delta}\Big\}\Big|_{x=0}.
\end{align*}

As a consequence of Proposition~\ref{InverseMellin} and the previous
discussion, we obtain
\begin{theorem}\label{ExpansionQ}
The trace of $Q(t)$ admits an asymptotic expansion
\begin{align*}
 \Tr Q(t) \sim_{t\to 0^+}\;
 & \sum_{k=0}^\infty a_k\, t^{(k-p)/m} +
 \sum_{k=0}^\infty \big\{ b_k + \beta_k \log t\big\} \,
 t^{k}\\
 & + \sum_{k=0}^\infty \big\{c_k + \gamma_k \log t + \delta_k (\log t)^2\big\}
 \, t^{z_k},
\end{align*}
where $z_k=\frac{k-m'-n}{m}$ with $n = \dim M$. Moreover,
\begin{align*}
 \beta_k &= -\frac{(-1)^k}{m k!}\rint{b} \res{\M(Q)(-k)}\\
 &\hspace{3em} - \frac{(-1)^k}{m k!(p+k m)!}\int_Y
 \partial_x^{p+k m}\Big\{x^{p + k m}\,
 \Res_0 \M(Q)(-k)|_{\Delta}\Big\}\Big|_{x=0}, \\
 \gamma_k &= -\frac{\Gamma_0(-z_k)}{m^2(p+z_k m)!}
 \int_Y \partial_x^{p+z_k m}\Big\{ x^{p +z_k m }
 \res{\M(Q)(-z_k)} \Big\}\Big|_{x=0}, \\
 \delta_k &= -\frac{(-1)^{z_k}}{m^2 z_k!(p+z_k m)!}
 \int_Y \partial_x^{p+z_k m}\Big\{x^{p + z_k m}\res{\M(Q)(-z_k)}
 \Big\}\Big|_{x=0}.
\end{align*}
If in any of the factorials $(p+k m)!$, $(p+z_k m)!$, or $(z_k)!$ the number
is not in $\N_0$, we define the corresponding coefficient to be $0$.
\end{theorem}

%%%%%%%%%%%%%%%%%%%%%%%%%%%%%%%%%%%%%%%%%%%%%%%%%%%%%%%%%%%%
\section{Heat trace expansion}\label{HeatTrace}
%%%%%%%%%%%%%%%%%%%%%%%%%%%%%%%%%%%%%%%%%%%%%%%%%%%%%%%%%%%%

Let $\Lambda$ be a sector of the form \eqref{sector}.
Let $P$ and $A$ be as in Theorem~\ref{GenResExpansion}. Thus,
\begin{equation}\label{splitting}
  P(A-\lambda)^{-1}=Q(\lambda)+R(\lambda)+S(\lambda),
\end{equation}
where the family $Q(\lambda)\in
 x^{m-p}\Psi_{c,\Lambda}^{m'-m,-m,m}(M,\Omega_b^{\frac 1{2}})$ is
holomorphically tempered, $R(\lambda)\in
 x^{-p}\Psi_{c,\Lambda}^{-\infty,\mathcal G(\alpha)}(M,\Omega_b^{\frac1{2}})$,
and $S(\lambda)\in
 x^{-p}\Psi_{\Lambda}^{-\infty,\mathcal F(\alpha)}(M,\Omega_b^{\frac 1{2}})$
for some index families $\mathcal{G}(\alpha)$ and $\mathcal
F(\alpha)$. Define $R(t)$ and $S(t)$ in the same way as $Q(t)$ was
defined in \eqref{IntegralQ}. Then by Theorems~3.7 and 3.8 of
\cite{LoHeatc01}, $R(t)$ and $S(t)$ are also operators of trace class and
we get a  decomposition
\[ \Tr(Pe^{-tA})=\Tr Q(t)+\Tr R(t)+\Tr S(t),\;\; t>0.\]

By the same theorems of loc.~cit. the trace of $R(t)$ admits an expansion
\[ \Tr R(t) \sim \sum_{k=0}^\infty r_k\, t^{(k-p)/m} \;\text{ as } t\to 0^+,\]
and $\Tr S(t)$ vanishes to infinite order at $t=0$. On the other
hand, the trace of $Q(t)$ admits the asymptotic expansion given in
Theorem~\ref{ExpansionQ}. Thus, $\Tr (Pe^{-tA})$ has the same
expansion as $\Tr Q(t)$. To provide a nicer expression for the
second term appearing in the formula for $\beta_k$ in
Theorem~\ref{ExpansionQ}, we proceed as follows. Denote by
$Q_0(\lambda)$ the function $Q(\lambda)$ in \eqref{splitting} for
$P = \mathrm{Id}$. Then, since $Q_0(\lambda)$ is equal to
$(A-\lambda)^{-1}$ modulo $\Psi_b^{-\infty}$, formally speaking,
\[ \M(Q_0)(z)=\frac i{2\pi}\int_\Upsilon \lambda^{-z} Q_0(\lambda)d\lambda
   \sim \frac i{2\pi}\int_\Upsilon \lambda^{-z}(A-\lambda)^{-1}d\lambda
   = A^{-z}. \]
Thus, although the complex powers $A^z$ does not exist in general,
we can still associate a useful meaning to it:
\begin{equation}\label{A^z}
   A^z := \M(Q_0)(-z).
\end{equation}
Moreover, the symbolic properties of $Q_0(\lambda)$ imply that
$\M(Q_0)(0)=\mathrm{Id}$ and that, modulo $\Psi_b^{-\infty}$,
$A^k \M(Q_0)(-z+k) = \M(Q_0)(-z)$ for any $k \in \N_0$.

Thus, we have proved:
\begin{theorem} \label{thm:main1}
Let $P\in x^{-p}\Psi_b^{m'}(M,\Omega_b^{\frac 1{2}})$, $p$,
$m'\in\R$. Let $A\in x^{-m}\Diff_b^m(M,\Omega_b^{\frac 1{2}})$ be
such that $A-\lambda$ is parameter-elliptic with respect to
$\alpha\in\R$ on a sector $\Lambda$ of the form \eqref{sector},
and assume that $m > p$. Then
\begin{align*}
 \Tr Pe^{-tA} \sim_{t\to 0^+}\;
 & \sum_{k=0}^\infty a_k\, t^{(k-p)/m} +
 \sum_{k=0}^\infty \big\{ b_k + \beta_k \log t\big\}\, t^{k}\\
 & + \sum_{k=0}^\infty \big\{c_k + \gamma_k \log t + \delta_k (\log t)^2\big\}
 \, t^{z_k},
\end{align*}
where $z_k=\frac{k-m'-n}{m}$ with $n = \dim M$. Moreover, the
coefficients $\beta_k$, $\gamma_k$ and $\delta_k$ are given explicitly by
\begin{align*}
\beta_k &= -\frac{(-1)^k}{m k!}\rint{b} \res{PA^k}\\
  &\hspace{3em} -\frac{(-1)^k}{m k!(p+k m)!}\int_Y
  \partial_x^{p+k m}\Big\{x^{p + k m}\,
  \Res_0 (PA^k A^z)|_\Delta|_{z=0} \Big\}\Big|_{x=0}, \\
\gamma_k &= -\frac{\Gamma_0(-z_k)}{m^2(p+z_k m)!}
  \int_Y \partial_x^{p+z_k m}\Big\{ x^{p +z_k m }
  \res{PA^{z_k}} \Big\}\Big|_{x=0}, \\
\delta_k &= -\frac{(-1)^{z_k}}{m^2 (z_k)!(p+z_k m)!}
  \int_Y \partial_x^{p+z_k m}\Big\{x^{p + z_k m}\res{P A^{z_k}}
  \Big\}\Big|_{x=0}.
\end{align*}
Again, if in any of the factorials the number is not in $\N_0$, we
define the corresponding coefficient to be $0$. The meaning of
the powers $A^{z}$ and $A^{z_k}$ is given in~\eqref{A^z}.
\end{theorem}

%%%%%%%%%%%%%%%%%%%%%%%%%%%%%%%%%%%%%%%%%%%%%%%%%%%%%%%%%%%%
\subsection*{Trace functionals revisited}
%%%%%%%%%%%%%%%%%%%%%%%%%%%%%%%%%%%%%%%%%%%%%%%%%%%%%%%%%%%%
By means of the generalized heat trace expansion obtained above,
we can recover the unique trace functionals on the algebras
$\mathcal{I}_{\sigma}$, $\mathcal{I}_{\partial}$, $\mathcal{A}_{\sigma}$,
$\mathcal{A}_{\partial}$ and $\mathcal{A}_{\sigma,\partial}$
from Section~\ref{TraceFunctionals}.

\begin{theorem} \label{thm:main2}
In the expansion of $\,\Tr Pe^{-tA}$ given in Theorem~\ref{thm:main1},
the coefficients $\beta_k$, $\gamma_k$ and $\delta_k$ can be written as
\begin{align*}
 \beta_k&= -\frac{(-1)^k}{m k!}
 \big( \Tr_\sigma(PA^k) + \Tr_\partial(PA^k)\big), \\
 \gamma_k &= -\frac{\Gamma_0(-z_k)}{m^2}
 \Tr_{\partial,\sigma}(PA^{z_k}), \\
 \delta_k &= -\frac{(-1)^{z_k}}{m^2 (z_k)!}
 \Tr_{\partial,\sigma}(PA^{z_k}).
\end{align*}
In particular, the coefficient of $\,\log t$ is
\begin{gather*}
  -\frac 1{m} \Tr_\sigma(P) -\frac 1{m}\Tr_\partial(P)
  -\frac 1{m^2} \Tr_{\partial,\sigma}(P), \\
\intertext{and}
  \Tr_{\partial,\sigma}(P)= -m^2\; \times
  \text{ the coefficient of $(\log t)^2$.}
\end{gather*}
\end{theorem}

\begin{corollary}
Suppose that $P\in x^{-p}\Psi^{m'}_b(M,\Omega_b^\frac12)$ with $p<0$. Then
there are no $(\log t)^2$ terms in the expansion of $\,\Tr Pe^{-tA}$, and
\[ \Tr_{\sigma}(P) = -m\;\times \text{ the coefficient of } \log t. \]
\end{corollary}

\begin{corollary}
Suppose that $P\in x^{-p}\Psi^{m'}_b(M,\Omega_b^\frac12)$ with $m'<-n$. Then
there are no $(\log t)^2$ terms in the expansion of $\,\Tr Pe^{-tA}$, and
\[ \Tr_{\partial}(P) = -m\;\times \text{ the coefficient of } \log t. \]
\end{corollary}

%% -------------------------------------------------------------
\bibliographystyle{amsplain}

\end{document}